\documentclass{article}

\usepackage[centertags]{amsmath} 
\usepackage{amssymb}
\usepackage{theorem,latexsym,epsfig,bbm,fancyheadings}

\input{xy}                       
\xyoption{all}

\theoremstyle{plain}                       
                                           
    \newtheorem{lemma}{Lemma}[section]  
    \newtheorem{kor}[lemma]{Corollary}     
    \newtheorem{satz}[lemma]{Theorem}      

\theorembodyfont{\rmfamily}                
                                          
    \newtheorem{Def}[lemma]{Definition}    
    \newtheorem{Bem}[lemma]{Remark}        

\newenvironment{beweis}{\noindent\textbf{Proof:}}{\par \hfill $\Box$ \hspace{1cm} \par}

\newenvironment{proofof}[1]{\renewcommand\proofname{#1} {\noindent\bf{Proof of #1:}}}%
{\par  \hfill {\bf{$\Box $}} \hspace{1cm}       \par}

\newcounter{zahl}%
\newenvironment{punkt}{\begin{list}{{\rm{(\roman{zahl})}}}%
    {\usecounter{zahl}%
     \setlength{\leftmargin}{0pt} \setlength{\itemindent}{4pt} \setlength{\topsep}{2pt} \setlength{\parsep}{2pt} }}%
    {\end{list}}%

\newcommand{\mc}[1]{\ensuremath{\mathcal{#1}}}
\newcommand{\sm}[2]{\ensuremath{#1\smash[b]{\wedge\,}#2}}
\newcommand{\smknapp}[2]{\ensuremath{#1\smash[b]{\wedge}#2}}%

\newcommand{\aus}{\raisebox{1pt}{\ensuremath{\,{\scriptstyle\in}\,}}}%
\newcommand{\eps}{\varepsilon}

\newcommand{\toh}[1]{\ensuremath{\stackrel{#1}{\rightarrow}}}%
\newcommand{\tow}[1]{\mbox{Turm-#1}}
\newcommand{\colim}{\operatorname*{colim}}
\newcommand{\holim}{\operatorname*{holim\,}}
\newcommand{\Rlim}{\operatorname*{lim^1}}
\newcommand{\frei}{\,\_\!\_\,}
\newcommand{\hofi}{\ensuremath{\mathfrak{Hofi}\,}}
\newcommand{\8}{\ensuremath{\infty}}                     
\newcommand{\Tot}{\operatorname{Tot}}

\newcommand{\kl}{\ensuremath{\raisebox{1pt}{\{}}}
\newcommand{\kr}{\ensuremath{\raisebox{1pt}{\}}}}
\newcommand{\hel}{\ensuremath{\raisebox{1pt}{[}}}
\newcommand{\her}{\ensuremath{\raisebox{1pt}{]}}}
\newcommand{\hrl}{\ensuremath{\raisebox{1pt}{(}}}
\newcommand{\hrr}{\ensuremath{\raisebox{1pt}{)}}}
\newcommand{\ho}[1]{\ensuremath{H\!o({#1})}}

\newcommand{\bth}{\raisebox{1pt}{\ensuremath{\,\scriptstyle \geq\,}}}%
\newcommand{\uber}{\raisebox{1pt}{\ensuremath{\,{\scriptstyle > }\,}}}%
\newcommand{\ul}[1]{\underline{#1}}
\newcommand{\ol}[1]{\overline{#1}}
\newcommand{\stabhom}[2]{\ensuremath{\hel#1,#2\her}}%
\newcommand{\Hom}[3]{\ensuremath{{\rm Hom}_{#1}\hrl#2,#3\hrr}}%
\newcommand{\gradhom}[4]{\ensuremath{{\rm Hom}_{#1}^{#2}\hrl#3,#4\hrr}}%
\newcommand{\Ext}[4]{\ensuremath{{\rm Ext}_{#1}^{#2}\hrl#3,#4\hrr}}%
\newcommand{\diag}[2]{ \begin{align} \begin{split} \xymatrix{#1} \end{split} \label{#2} \end{align}}%

\numberwithin{equation}{section}        
\newcommand{\Ref}[1]{\hrl\ref{#1}\hrr}        

\newcommand{\nummer}[2]{ \begin{equation} \begin{split} #1 \end{split} \label{#2} \end{equation} }   

\def\Sig{\Sigma}
\def\lam{\lambda}		
\def\kap{\kappa}
\def\gam{\gamma}

\author{Georg Biedermann}
\title{Injective completion with respect to homology}

\begin{document}

\maketitle

\begin{abstract}
Generalizing $F$-nilpotent completion for a ring spectrum $F$ we first define the notion of completion with respect to a thick subcategory in a monogenic stable homotopy category.
Specializing this to the thick subcategory generated by $F$-injectives gives an injective completion functor. This is the completion functor adapted to the modified Adams spectral sequence, which uses absolute instead of relative injective resolutions. Finally we show, that both constructions coincide for suitable ring spectra.
\end{abstract}
\vspace{1cm}

{\bf\Large Introduction} 
\vspace*{2ex}

Completions in algebraic topology can be used to study the convergence behaviour of spectral sequences. 
In order to investigate the convergence of the Adams spectral sequence Bousfield in \cite{Bou:loc} introduced the notion of nilpotent completion of a spectrum with respect to a generalized homology theory. This notion was especially adapted to the original construction by Adams. There are other constructions of Adams spectral sequences. 
The construction via injective resolutions was introduced in \cite{Bri:adams}, and more recently considered in \cite{Bou:klocal} and \cite{Fra:uni}. It is called the modified Adams spectral sequence. 
The reader is also referred to \cite{Dev:brown-comenetz} for more results and applications.

Analogously to nilpotent completion I define the notion of injective completion which is adapted to the construction of the modified Adams spectral sequence using injective resolutions. 
To give a unifying exposition I found it very convenient to use the language of monogenic stable homotopy categories in the sense of \cite{HPS:Ax}. One can formulate the concept of completion with respect to a thick subcategory in this frame work. Nilpotent and injective completion are then just special cases where we pick the appropriate thick subcategory.

To be able to compare these two completions I use the theory of homotopical completion from \cite{Bou:cos}. I have included a short summary of the relevant definitions and theorems from there.

The main result is \ref{inj=nil}. It asserts that injective completion and nilpotent completion are equivalent. I hope this theorem gives information on convergence questions for the modified Adams spectral sequence based on ring spectra satisfying our assumptions.

The theorem \ref{inj=nil} is a corollary of theorem \ref{bousfield} of Bousfield which tells us that the completeness result \ref{resultat} is a sufficient condition for \ref{inj=nil}. 

The first section is a collection of general statements mostly needed to fix notation. 
The next two sections are devoted to the construction of the modified Adams spectral sequence. Section 4 and 5 outline the two view points of general stable homotopical completion. In section 6 I define nilpotent and injective completion. The main technical part and the preparation for the later proofs are put in section 7 and 8. The main theorem is in section 9 and some immediate applications are in section 10.

This article is a more developped version of my Diplomarbeit presented to the Rheinische Friedrich-Wilhelms-Universit\"at Bonn in April 2000. I thank my advisor Jens Franke for drawing my attention to this problem and for all the good comments at the right time. I thank E.S. Devinatz for pointing out to me an error in a previous version. I am grateful to A.K. Bousfield for explaining me the ideas in his preprint \cite{Bou:cos} and how they are connected to what I did. 


\tableofcontents \pagebreak

\section{Monoids in stable homotopy categories}

Let \mc{D} be a stable homotopy category in the sense of \cite{HPS:Ax} and let  $F$ be a monoid in this category. The following facts taken from \cite{Ada:lec} or \cite{Rav:greenbook} immediately generalize to this situation.

\begin{Def}
\begin{sloppypar}
Let {\bf\boldmath $\pi_*F$-mod\unboldmath} be the category of $\pi_*F$-left modules.
$\Hom{\pi_*F}{\frei}{\frei}$ shall denote the morphisms in this category. 
\end{sloppypar}
\end{Def}

\begin{lemma} 
\label{Künneth}
Let $F$ be a monoid in \mc{D}. Assume that $F_*F$ is flat as a right module over $\pi_*F$. Consider the map
\begin{center}
      $\xymatrix{ \sm{(\sm{F}{F})}{(\sm{F}{X})} \ar[rr]^{\smknapp{1}{\smknapp{\mu_F}{1}}} & & \sm{F}{\sm{F}{X}}  .} $
\end{center}
It induces a map 
\begin{center}
    $F_*F \otimes_{\pi_*F} F_*X \to F_*\sm{F}{X} .$
\end{center} 
With our assumptions it is a natural isomorphism of $\pi_*F$-modules.
\end{lemma}

\begin{Def}
We denote by {\bf\boldmath $F_*F$-comod\unboldmath} the category of $F_*F$-left comodules. We denote the homomorphisms and Ext-groups by $\Hom{F_*F}{\frei}{\frei}$ and $\Ext{F_*F}{}{\frei}{\frei}$.
\end{Def}

\begin{Bem}
We know that if $F_*F$ is flat as a right module over $\pi_*F$, then $(\pi_*F,F_*F)$ is a Hopf algebroid, the category $F_*F$-comodules is abelian, the forgetful functor $F_*F$-comod $\to \pi_*F$-mod possesses a right adjoint, which is given by $F_*F \otimes_{\pi_*F} \frei $ and that the functor $F_*$ takes values in the category $F_*F$-comod. 
It also follows that $F_*F$ is faithfully flat as $\pi_*F$-right module.
\end{Bem}

\begin{Def}  
\label{Def. univ. Koeff.}
Let $M$ be an $F$-module. We assume that $F$ is the colimit of finite objects $F_\alpha$, satisfying the following properties:
\begin{punkt}
    \item $F_*(DF_\alpha)$ is projective over $\pi_*F$.
    \item For every $\alpha$ the map
       \begin{center}
          $ F_*F_\alpha \otimes_{\pi_*F} \pi_*M \to M_*F_\alpha \, , $
       \end{center}
          which is induced by the composition $\xymatrix{\sm{F}{\sm{F_\alpha}{M}} \ar[r]^{\cong} & \sm{F}{\sm{M}{F_\alpha}} \ar[r]^{\sm{\mu_M}{1}} & \sm{M}{F_\alpha}}$ is an isomorphism.
\end{punkt}
We say that $F$ has the {\bf universal coefficient property}, if for every $F$-module $M$ these properties are satisfied. Recently some authors call this property topologically flat.
\end{Def}

\begin{satz} 
\label{univ. Koeff.}
Let $F$ satisfy the universal coefficient property. Then there is for all $X$ in \mc{D} a spectral sequence
\begin{center}
    $ E_2^* = F_*X \otimes_{\pi_*F} \pi_*M\,\, \Longrightarrow\,\, M_*X$
\end{center}
If $F_*X$ is flat, then the edge homomorphism
\begin{equation*}
      F_*X \otimes_{\pi_*F} \pi_*M \to M_*X \label{Kante}
\end{equation*}
is an isomorphism of $\pi_*F$-modules. 
\end{satz}

\begin{kor} \label{m}
The map
\begin{equation}
    F_*F \otimes_{\pi_*F} \pi_*M \stackrel{\cong}{\longrightarrow} F_*M , \notag
\end{equation}
is an isomorphism for every $F$-module $M$.
\end{kor}

\begin{Def} \label{Dualität}
Let $F$ be the colimit of finite objects $F_\alpha$ satisfying the following properties:
\begin{punkt}
    \item $F_*DF_\alpha$ is finitely generated and projective over $\pi_*F$.
    \item For every $\alpha$ and every $F$-module $M$ the homomorphism
       \begin{center}
           $ M^*DF_\alpha \to \Hom{\pi_*F}{F_*DF_\alpha}{\pi_*M} $,
       \end{center} 
mapping $DF_\alpha \to M$ to $F_*DF_\alpha \to F_*M \stackrel{(\mu_M)_*}{\longrightarrow} \pi_*M$ is an isomorphism of $\pi_*F$-modules. 
\end{punkt}
We will say, that $F$ has the {\bf duality property}, if $(i)$ and $(ii)$ are satisfied.
\end{Def}

\begin{satz} 
\label{univ. Dual.}
Let $F$ satisfy the duality property. Let $F_*X$ be projective over $\pi_*F$. Then there is for every $F$-module $M$ an isomorphism
\begin{equation*}
    M^*X \stackrel{\cong}{\longrightarrow} \emph{\gradhom{\pi_*F}{*}{F_*X}{\pi_*M}}.
\end{equation*}
In particular we obtain an isomorphism $F_*DF_{\alpha}\cong\emph{\Hom{\pi_*F}{F_*F_{\alpha}}{\pi_*F}}$.
\end{satz} 

\begin{Bem}
The duality property is property 1.1 from \cite{Dev:brown-comenetz}. The universal coefficient property is essential in the proof of \ref{Behauptung 2}. Both properties are satisfied if $F$ is an evenly graded and Landweber exact ring spectrum as shown in \cite{Dev:brown-comenetz}.
\end{Bem}

\section{Injective objects and Eilenberg-MacLane objects}

\begin{Def}  
\label{Eilenberg-MacLane}
Let $F_*$:$\mc{D} \to \mc{A}$ be a homological functor from a triangulated category to an abelian category with a canonical equivalence $F_*\Sigma X\cong F_{*-1}X$, where $\Sigma$ is the shift functor in \mc{D}. Consider the following functor
\begin{equation}
      X \mapsto \Hom{\mc{A}}{F_*X}{I} \nonumber
\end{equation} 
If this functor is representable by an object $E_I$ from \mc{D} and if the canonical morphism $F_*E_I\to I$ induced by
\begin{equation}
      id_{E_I}\in\Hom{\mc{D}}{E_I}{E_I} \cong \Hom{\mc{A}}{F_*E_I}{I} \nonumber
\end{equation}
is an isomorphism, then we call $E_I$ an {\bf\boldmath $(F,I)$-Eilenberg-MacLane object\unboldmath}. Sometimes we will write $E(I)$.

Let $I$ be some injective object in \mc{A}. Then we call an $(F,I)$-Eilenberg-MacLane object an {\bf\boldmath $F$-injective object\unboldmath}. 

By abuse of language we will often make no difference between corresponding elements in the isomorphic groups $\Hom{\mc{D}}{X}{E_I}$ and $\Hom{\mc{A}}{F_*X}{I}$, although they are morphisms in different categories.
\end{Def}

\begin{Def} 

We will say that the functor $F_*$ is {\bf injectively resolvable}, if every object in \mc{D} admits a morphism to an $F$-injective object, that induces a monomorphism in \mc{A}.
\end{Def}

\begin{lemma}  \label{spez-inj-Spek-F-Modul}
Let \mc{D} be a monogenic stable homotopy category.
Let $F$ be a monoid having the duality property. Let $I$ be an injective $\pi_*F$-left module, and let $E_I$ be an object representing the functor 
\begin{equation}
    X \mapsto \emph{\Hom{\pi_*F}{F_*X}{I}} \cong \emph{\Hom{F_*F}{F_*X}{F_*F \otimes_{\pi_*F} I}}    \nonumber
\end{equation}
Then $E_I$ is an $F$-module, and there is an isomorphism 
\begin{equation*}
    F_*E_I \cong F_*F \otimes_{\pi_*F} I  
\end{equation*}
of $F_*F$-comodules.
In particular is $E_I$ an $(F_*,F_*F \otimes_{\pi_*F} I)$-Eilenberg-MacLane object and hence $F$-injective for the category of $F_*F$-comodules. If $F$ also satisfies the universal coefficient property, then we have
\begin{equation*}
    \pi_*E_I \cong I  
\end{equation*}
as $\pi_*F$-modules.
\end{lemma}

\begin{beweis}
For the proof of this fact the reader is referred to \cite{Fra:uni} or to \cite{Dev:brown-comenetz}.
\end{beweis}

\begin{kor} \label{Tensor-Auflösung}
Every $F_*F$-comodule admits a monomorphism to an injective $F_*F$-co\-mo\-dule, for which there exists an Eilenberg-MacLane object. These objects can be chosen of the form $F_*F \otimes_{\pi_*F} I$, where $I$ is an injective $\pi_*F$-module. In particular $F_*:\mc{D} \to F_*F$-comod is injectively resolvable.
\end{kor}

\section{Construction of an injective resolution}

Let $F_*:\mc{D}\to\mc{A}$ be an injectively resolvable functor.
Let $Y$ be in $\mc{D}$. Then there is an injective resolution  
\nummer{  0 \to F_*Y \toh{(\eps^0)_*} I^0 \toh{i^0} I^1 \toh{i^1} \ldots }{AuflF_*Y}
in \mc{A}, such that for every $I^s$ there exists an Eilenberg-MacLane object $E(I^s)$. One constructs an injective resolution of $Y$ in \mc{D} with the following inductive machinery:  
Set $Y^0:=Y$.
After having constructed $Y^s$ together with an injection \mbox{$ 0 \to F_*Y^s \toh{(\eps^s)_*} I^s $}, this induces a map \mbox{$ Y^s \toh{\eps^s} I^s $} in \mc{D}. We define $Y^{s+1}$ as the fiber of $\eps^s$ and form the triangle 
\nummer{  Y^{s+1} \toh{\gam^s} Y^s \toh{\eps^s} E(I^s) \to \Sig Y^{s+1} . }{Dreieps^s}
Obviously $F_*\gam^s = 0$ holds, so that the long exact $F_*$-sequence splits into the following short exact sequences:
\nummer{  0 \to F_*Y^s \toh{(\eps^s)_*} I^s \to F_{*-1}Y^{s+1} \to 0 }{F-Seqeps^s}
The map $I^s \toh{i^s} I^{s+1} $ factors as follows:
\diag{ & & & 0 \ar[d] & \\
       0 \ar[r] & F_*Y^s \ar[r]^{(\eps^s)_*} & I^s \ar[r] \ar[dr]_{i^s} & F_{*-1}Y^{s+1} \ar[r] \ar[d]^{(\eps^{s+1})_*} & 0 \\              & & & \Sigma I^{s+1} & }{i^sfaktor} 
The construction is finished producing the following diagram:
\diag{ Y \ar@{=}[r] & Y^0 \ar[d]_{\eps^0} & Y^1 \ar[l]_{\gam^0} \ar[d]_{\eps^1} & Y^2 \ar[l]_{\gam^1} \ar[d]_{\eps^2} & \ldots \ar[l]_{\gam^2} \\
       & E(I^0) \ar[ur]_+ & E(I^1) \ar[ur]_+ & E(I^2)\ar[ur]_+  & }{AuflüberY}
where $F_*\gam^s = 0$ for all $s\bth 0$. Here $+$ means as in the sequel, that the map has in fact degree $-1$.

\begin{Def}
The diagram \Ref{AuflüberY} is called an {\bf\boldmath $F$-injective resolution\unboldmath} or shortly injective resolution {\bf over $Y$}.
\end{Def}

From this resolution over $Y$ we derive a resolution under $Y$. This procedure is not depending on $F$-injectivity. So let $\{Y^s\}$ be an arbitrary resolution over $Y$.

Set $Y_{-1} = 0$, and for $s \bth 1$ define $Y_{s-1}$ by
\nummer{  Y^s \to Y^0 \to Y_{s-1} \to \Sig Y^s , }{DreiDefY_s-1}
where $ Y^s \to Y^0 $ is the composition of the tower maps $ \gam^s $. To construct a map $\xymatrix@1{Y_s \ar[r]^-{\gam_{s-1}} & Y_{s-1}}$, we use the octahedron axiom for the composition $Y^{s+1} \to Y^s \to Y^0 $. We obtain a diagram
\diag{ & \ar[d]_+ & \ar[d]_+ &  \\
      Y^{s+1} \ar[r]^{\gamma^s} \ar@{=}[d] & Y^s \ar[r]^{\eps^s} \ar[d] & E(I^s) \ar[r]^+ \ar[d] & \\
      Y^{s+1} \ar[r] & Y^0 \ar[r]  \ar[d] & Y_s \ar[r]^+ \ar[d]^{\gamma_{s-1}} & \\
      & Y_{s-1} \ar@{=}[r] & Y_{s-1} & }{Oktüberunter}
in which each column and row is a distinguished triangle together with the morphisms of triangles
\diag{ Y^{s+1} \ar[r] \ar[d]^{\gam^s} & Y \ar[r] \ar@{=}[d] & Y_s \ar[r] \ar[d]^{\gam_{s-1}} & \Sig Y^{s+1} \ar[d]^{\Sig \gam^s} \\
       Y^s \ar[r] & Y \ar[r] & Y_{s-1} \ar[r] & \Sig Y^s }{überunter}
and
\diag{ E(I^s) \ar[r] \ar@{=}[d] & Y_s \ar[r]^{\gam_{s-1}} \ar[d] & Y_{s-1} \ar[r] \ar[d] & \Sig E(I^s) \ar@{=}[d] \\
       E(I^s) \ar[r] & \Sig Y^{s+1} \ar[r]^{\Sig \gam^s} & \Sig Y^s \ar[r] & \Sig E(I^s)\ .}{Morüberunter} 
Altogether we arrive at a diagram
\diag{ Y \ar@{=}[r] \ar[d] & Y \ar@{=}[r] \ar[d] & Y \ar@{=}[r] \ar[d] & \ldots \\
           Y_0 & Y_1 \ar[l]_{\gam_0} & Y_2 \ar[l]_{\gam_1} & \ldots \ar[l]_{\gam_2}}{AuflunterY}
 
\begin{Def}  \label{associated}
We call this the {\bf resolution under $Y$ associated to $\{Y^s\}$}. Note that the maps $Y \to Y_s$ are considered as part of the structure (see Def. \ref{pro}). Note also that this process works equally the other way around, this is then called {\bf the associated tower over $Y$}.  
\end{Def}

\begin{Def}  \label{injektive Auflösung}
If $\{Y^s\}$ was an injective resolution over $Y$, then the  associated resolution is called an {\bf\boldmath $F$-injective resolution\unboldmath} or for short injective resolution {\bf under $Y$}. 
\end{Def}

\begin{Def}  \label{modifizierte ASS}
We rewrite the diagram \Ref{Morüberunter}, so an injective resolution under $Y$ is a sequence of distinguished triangles of the form:
\diag{ \ldots & Y_{s-1} \ar[l] \ar[d]_+ & Y_s \ar[l]_{\gam_{s-1}} \ar[d]_+ & Y_{s+1} \ar[l]_{\gam_{s}} & \ldots \ar[l] \\ 
        & E(I^s) \ar[ur] & E(I^{s+1}) \ar[ur] & }{liefert exaktes Paar}
If we apply the functor $\stabhom{X}{\frei}$ to this diagram, we get an unraveled exact couple. The spectral sequence derived from it is called the {\bf modified Adams spectral sequence}.
\end{Def}


\section{Bousfield's \mc{G}-completion}

The original theory of \cite{Bou:cos} is here adjusted to a stable situation (see Def. \ref{linear}) by substituting $n\bth 0$ by $n\aus\mathbbm{Z}$ and leaving out all loop functors in all conditions .

\begin{Def}
Let \mc{C} be a simplicial left proper stable model category, and let \mc{G} be a class of objects in the homotopy category \ho{\mc{C}}. (Originally these objects were required to be group objects, but since we are in the stable case every object is an (abelian) group object.) A map $i:A\to B$ in \ho{\mc{C}} is called {\bf \mc{G}-monic} when $i^*:\stabhom{B}{G}_n\to\stabhom{A}{G}_n$ is surjective for each $G\aus\mc{G}$ and each $n\aus\mathbbm{Z}$.

An object $I$ is called {\bf \mc{G}-injective} when $i^*:\stabhom{B}{I}_n\to\stabhom{A}{I}_n$ is surjective for each \mc{G}-monic map $i:A\to B$ and each $n\aus\mathbbm{Z}$.

We say that \ho{\mc{C}} {\bf has enough \mc{G}-injectives} if each object in \ho{\mc{C}} is the source of a \mc{G}-monic map to a \mc{G}-injective target. We say that \mc{G} is {\bf functorial}, if these maps can be chosen functorially.
\end{Def}

\begin{Def}
Let $c\mc{C}$ be the category of cosimplicial objects over \mc{C}. We call a map $f:X^{\bullet}\to Y^{\bullet}$ a
\begin{punkt}
    \item {\bf \mc{G}-equivalence} if $f_*:\stabhom{Y^{\bullet}}{G}_n\to\stabhom{X^{\bullet}}{G}_n$ is a weak equivalence of simplicial groups for each $n\aus\mathbbm{Z}$.
    \item {\bf \mc{G}-cofibration} if $f$ is a Reedy-cofibration and $f^*:\stabhom{Y^{\bullet}}{G}_n\to\stabhom{X^{\bullet}}{G}_n$ is a fibration of simplicial groups for each $G\aus\mc{G}$ and $n\aus\mathbbm{Z}$.
    \item {\bf \mc{G}-fibration} if $f:X^n\to Y^n\times_{M^nY^{\bullet}}M^nX^{\bullet}$ is a \mc{G}-injective fibration for $n\aus\mathbbm{Z}$.
\end{punkt}
\end{Def}

\begin{satz}
The category $c\mc{C}$ of cosimplicial objects over a simplicial left proper stable model category \mc{C} with a class \mc{G} and enough \mc{G}-injectives together with the previously described classes of maps becomes a simplicial left proper pointed model category.
\end{satz}

\begin{Def}
A {\bf \mc{G}-resolution} of an object $A$ in \mc{C} consists of an acyclic \mc{G}-cofibration $A\to\bar{A}^{\bullet}$ to a \mc{G}-fibrant object $\bar{A}^{\bullet}$ in $c\mc{C}$ where $A$ is considered constant in $c\mc{C}$.
\end{Def}

\begin{Def}
For an object $A$ in \mc{C} we define the {\bf \mc{G}-completion} $\alpha:A\to\widehat{L}_{\mc{G}}A$ in \ho{\mc{C}} by setting
\begin{equation*}
    \widehat{L}_{\mc{G}}A:= {\rm Tot} \bar{A}^{\bullet}
\end{equation*}
where $A\to\bar{A}^{\bullet}$ is a \mc{G}-resolution of $A$. 
\end{Def}

\begin{Bem}
This determines a functor $\widehat{L}_{\mc{G}}:\mc{C}\to\ho{\mc{C}}$ which is well-defined up to natural equivalence. In fact the \mc{G}-completion will give a functor $\widehat{L}_{\mc{G}}:\ho{\mc{C}}\to\ho{\mc{C}}$ and a natural transformation $\alpha: Id\to\widehat{L}_{\mc{G}}$ belonging to a monad on \ho{\mc{C}}. When \mc{C} has functorial replacements and \mc{G} is functorial, the \mc{G}-completion is canonically represented by a functor $\widehat{L}_{\mc{G}}:\mc{C}\to\mc{C}$.
\end{Bem}

\begin{satz}
Suppose $A\to Y^{\bullet}$ is a weak \mc{G}-resolution and $\ul{Y}^{\bullet}$ is a Reedy-fibrant replacement of $Y^{\bullet}$. Then there is a natural equivalence
\begin{equation*}
    \widehat{L}_{\mc{G}}A\cong {\rm Tot}\ul{Y}^{\bullet}
\end{equation*}
in \ho{\mc{C}}.
\end{satz}

\begin{Def}
An object $A$ in \ho{\mc{C}} is called {\bf \mc{G}-complete} if $\alpha:A\to\widehat{L}_{\mc{G}}A$ is an isomorphism in \ho{\mc{C}}.

A {\bf \mc{G}-complete expansion} of an object $A$ in \mc{C} consists of a \mc{G}-equivalence $A\to Y^{\bullet}$ in $c\mc{C}$ such that $Y^n$ is \mc{G}-complete for $n\bth 0$.
\end{Def}

\begin{satz}
If $A\to Y^{\bullet}$ is a \mc{G}-complete expansion and $\ul{Y}^{\bullet}$ is a functorial Reedy-fibrant replacement of $Y^{\bullet}$, then there is a natural equivalence
\begin{equation*}
    \widehat{L}_{\mc{G}}A\cong {\rm Tot}\ul{Y}^{\bullet}
\end{equation*}
in \ho{\mc{C}}.
\end{satz}

\begin{satz} \label{bousfield}
Suppose \mc{G} and $\mc{G}'$ are classes of injective models in $\ho{\mc{C}}$. If each \mc{G}-injective object is $\mc{G}'$-injective and each $\mc{G}'$-injective object is \mc{G}-complete, then there is a natural equivalence
\begin{equation*}
    \widehat{L}_{\mc{G}}X\cong\widehat{L}_{\mc{G}'}X
\end{equation*}
for every $X$ in \ho{\mc{C}}.
\end{satz}

\section{Completion with respect to thick sub\-ca\-te\-gories}

This paragraph is just the result of observing that completion as in \cite[par. 5]{Bou:loc} works in the case of monogenic stable homotopy categories in the sense of \cite{HPS:Ax}.

\begin{Def}  
\label{pro}
Let \mc{K} be a category. We define the category \tow{\mc{K}} as follows:
\begin{punkt}
      \item Objects are towers in \mc{K}, which are diagrams of the form
         \begin{equation}
              { Y^0 \leftarrow Y^1 \leftarrow Y^2 \leftarrow \ldots } \nonumber
         \end{equation}
      \item  $\Hom{}{\{X_s\}}{\{Y_s\}} := \lim\limits_{t}\, \colim\limits_{s}\, \Hom{\mc{K}}{X_s}{Y_t}$        
\end{punkt}
There is a canonical functor
\begin{center}
 $ \begin{array}{cc}
       \mc{K} \to  & \tow{\mc{K}}  \\
          Y \mapsto  & \{Y\}      
   \end{array}                      $
\end{center}   
taking $Y$ to the constant tower $ Y = Y = Y = \ldots $ denoted by $\{Y\}$. In this way \mc{K} becomes a full subcategory of \tow{\mc{K}}. We call a tower $\{Y^s\}$ a {\bf tower over $Y$}, if $Y=Y^0$. We call $\{Y_s\}$ a {\bf tower under $Y$}, if there is a map $\{Y\}~\to~\{Y_s\}$ in \tow{\mc{K}}. We will write the indices of towers in correspondence with their property of being a tower over or under $Y$ with the exception of lemma \ref{Tot} and its proof, where also cosimplicial objects are considered. 

\begin{sloppypar}	
Morphisms of towers, that are isomorphisms in this category, are called {\bf pro-isomorphisms}. This is the terminology of \cite{BK:lim}. 
\end{sloppypar}
\end{Def}

\begin{Def}  
\label{B-Aufl}
Let \mc{B} be a thick subcategory of a triangulated category \mc{D}. A tower $\{W_s\}$ under an object $W$ is called a {\bf \mc{B}-resolution} of $W$, if the following properties are satisfied:
\begin{punkt}
    \item For every $s \bth 0$ $W_s$ is an object in \mc{B}. 
    \item For every $N$ in \mc{B} the canonical map $\colim_{s} \stabhom{W_s}{N} \to \stabhom{W}{N}$ is an isomorphism.
\end{punkt}      
\end{Def}

\begin{Def}  
We say, that a full subcategory \mc{B} is {\bf small up to isomorphism} with representatives \mc{V}, if there exists a (small) set \mc{V} of objects of \mc{B}, such that every object in \mc{B} is isomorphic to an object of \mc{V}.
\end{Def}

\mc{B}-resolutions may not exist, however (inspired by \cite{Bou:cos} and \cite{HPS:Ax}):
\begin{lemma} 
If \mc{B} is small up to isomorphism and colocalizing, meaning that it is closed under arbitrary products, then every object has a \mc{B}-resolution.
\end{lemma}

\begin{beweis}
Let \mc{V} be a set of representatives of \mc{B}, then there is an obvious map 
\begin{equation*}
      X\to\prod_{B\in\mc{V}}\,\,\prod_{n\in\mathbb{Z}}\,\,\prod_{f\in[X,B]_n}\!\! B\,\, =:\, B_0=:X_0. 
\end{equation*}
Take the fiber to form a distinguished triangle and repeat the construction. We get a tower over $X$, and by applying the upside-down-construction of \Ref{DreiDefY_s-1} we obtain a tower $\{X_s\}$, that is the associated tower under $X$ (Def. \ref{associated}).
Inductively we see, that all $X_s$ are in \mc{B}. Because there is a long exact sequence
\begin{equation*}
    ...\to\colim_s\stabhom{X^s}{B}_{*+1}\to\colim_s\stabhom{X_s}{B}_*\to\stabhom{X}{B}_*\to\colim_s\stabhom{X^s}{B}_*\to...
\end{equation*}
it suffices to show $\colim_s\stabhom{X^s}{B}=0$ for every $B$ in \mc{B} to prove the second condition. But this is obvious, since the map $\stabhom{X^s}{B}\to\stabhom{X^{s+1}}{B}$ vanishes. So $\{X_s\}$ is a \mc{B}-resolution.
\end{beweis}

\begin{lemma} \label{pi-Iso}
Let \mc{B} be a thick subcategory of \mc{D}. Let $\{W_s\}$ and $\{V_s\}$ be \mc{B}-resolutions of $Y$. Then there exists a unique pro-isomorphism $e:\{W_s\} \to \{V_s\}$ in \tow{\mc{D}}, such that
\begin{center}
      $\xymatrix{ \{Y\} \ar@{=}[r] \ar[d] & \{Y\} \ar[d] \\
                  \{W_s\} \ar[r]^e & \{V_s\} }$
\end{center}
commutes.
\end{lemma}

\begin{Def}  

If \mc{D} is a monogenic stable homotopy category (see \cite{HPS:Ax}), then $\stabhom{S}{\frei}_* =: \pi_*$ is a homological functor to abelian groups with the property that $\pi_*(X) = 0   \Rightarrow  X \cong 0$. A map $g:\{W_s\} \to \{V_s\}$ from \tow{\mc{D}} is called {\bf weak equivalence}, if the induced map $\pi_ig: \{\pi_iV_s\} \to \{\pi_iW_s\}$ is a pro-isomorphism in \tow{Ab} for every $i \aus \mathbb{Z}$. 

\end{Def}

The next lemma is \cite[Lemma 5.11]{Bou:loc}.
\begin{lemma} 
\label{Iso der Homotopielimites}
A morphism in \tow{\mc{D}} induces a map of homotopy limits.
Let $g:\{V_s\}\!\to\!\{W_s\}$ be a weak equivalence and let $V_{\8}, W_{\8}$ be the homotopy limits of two towers. There is a map \mbox{$u:V_{\8} \to W_{\8}$}, such that
\begin{center}
    $\xymatrix{ V_{\8} \ar[r]^{\{u\}} \ar[d] & W_{\8} \ar[d] \\
                \{V_s\} \ar[r]^g  &  \{W_s\}  }    $
\end{center} 
commutes in \tow{\mc{D}}, and this $u$ is an isomorphism in \mc{D}.
\end{lemma}

\begin{kor} 
\label{dicke Vervollst} 
Let \mc{D} be a monogenic stable homotopy category, and let $Y$ be an arbitrary object in \mc{D}. If $\{W_s\}$ and $\{V_s\}$ are two \mc{B}-resolutions of $Y$, then:
\begin{align} 
      \holim_s W_s = \holim_s V_s \nonumber
\end{align}
\end{kor}

\begin{beweis}
By lemma \ref{pi-Iso} there is an isomorphism $\{W_s\} \to \{V_s\}$ in \tow{\mc{D}}, this is in particular a weak equivalence. Now the assertion follows from lemma \ref{Iso der Homotopielimites}, where we have shown that weak equivalences induce isomorphisms of homotopy limits.
\end{beweis}

\begin{Def}  

In view of the last lemma we write
\begin{align}
       \holim_s W_s =: \mc{B}^{\wedge}W \notag
\end{align}
for a thick subcategory \mc{B} of a monogenic stable homotopy category \mc{D} and for a \mc{B}-resolution $\{W_s\}$ of an object $W$ in \mc{D}. We call $\mc{B}^{\wedge}W $ a {\bf \mc{B}-completion} of $W$. Clearly \mc{B}-completion commutes with $\Sig^{n}$. It exists if at least one \mc{B}-resolution for $W$  exists. 
\end{Def} 

\begin{Bem} \label{triangulated case}
Here we can use the definition of $\holim$ of \cite{BN:lim} or that of a sequential limit of \cite{HPS:Ax}.
As long as we stay in the triangulated context there will always be difficulties with functoriality here. In this section and in the next one we will indicate this problem in the statements. At the price of choosing an underlying model category as in the later sections it can be overcome.
\end{Bem}

\begin{lemma} \label{Dualität der Vervollst.}
There is a map $W \to \mc{B}^{\wedge}W$ which is in general not natural. Let $\{W^s\}$ be a resolution over $W$, such that its associated resolution under $W$ is a \mc{B}-resolution. Then there are \emph{(}generally not natural\emph{)} maps
\begin{equation}
    \holim_s W^s \to W \to \mc{B}^{\wedge}W \to \Sig \holim_s W^s \, ,\nonumber
\end{equation}
forming a distinguished triangle. The isomorphism class of $\holim W^s$ is independent of a choice of a resolution over $W$.
\end{lemma}

\begin{Def} \label{dualer Term}
Let $\{W^s\}$ be a resolution over $W$, such that its associated resolution is a \mc{B}-resolution. We call
\begin{equation}
    W^{\wedge}\mc{B} := \holim_s W^s \nonumber
\end{equation}
the {\bf opposite term}. Lemma \ref{Dualität der Vervollst.} shows that this notion is well defined. 
\end{Def}

\begin{lemma} \label{Eigenschaften} 
\mc{B}-completion and its opposite term are exact in the sense that there is a distinguished triangle of completions if there was initially one. They preserve finite coproducts and retracts.
\end{lemma}

The relation between this concept of completion and the one in section 4 is given by

\begin{lemma} \label{Tot}
Let \mc{G} be a class of injective models in a left proper simplicial stable model category \mc{C} such that \ho{\mc{C}} is monogenic. Let $\mc{B}(\mc{G}{\rm -inj})$ be the thick subcategory of $\ho{\mc{C}}$ generated by the \mc{G}-injective objects. Let $Y\to Y^\bullet$ be a fibrant replacement of an object $Y$ of \mc{C} in the \mc{G}-model structure on $c\mc{C}$. Then
\begin{equation*}
    \Tot_0Y^\bullet\leftarrow ... \leftarrow\Tot_sY^\bullet\leftarrow\Tot_{s+1}Y^\bullet\leftarrow ...
\end{equation*} 
is a $\mc{B}(\mc{G}{\rm -inj})$-resolution under $Y$. Therefore:
\begin{equation*}
    \widehat{\mc{L}}_{\mc{B}(\mc{G}{\rm -inj})}Y=\Tot Y^\bullet\cong\mc{B}(\mc{G}{\rm -inj})^{\wedge}Y  
\end{equation*}
\end{lemma}

The following proof is due to Bousfield.

\begin{beweis}
We have to prove the two conditions from Def. \ref{B-Aufl}. $(i)$ follows by induction from the definition of $\Tot_s$. For $(ii)$ let $I$ be a \mc{G}-injective object and consider the $I$-cohomology spectral sequence of the tower $\{\Tot_sY^\bullet\}$.
Its $E_1$-term consists of
\begin{equation*}
    E_1^s=I^*Fib_s=\stabhom{Fib_s}{I}_*,
\end{equation*}
where $Fib_s$ is the fiber of $\Tot_s Y^\bullet\to\Tot_{s-1}Y^\bullet$. By \cite[p. 282]{BK:lim} or \cite[p. 391]{GoJar:simp} there is an isomorphism $Fib_s\cong \Omega^sN^sY^\bullet$, where $N^sY^\bullet:=\text{fiber}(Y^s\to M^{s-1}Y^\bullet)$ is the geometric normalization of $Y^\bullet$. 
Moreover it is true, that there is an isomorphism
\begin{equation*}
    I^*(Fib_s)=I^*(\Omega^sN^sY^\bullet)\cong N^s(I^{*-s}Y^\bullet),
\end{equation*}
where on the right hand side $N^s$ denotes the usual normalization of complexes.
Also the spectral sequence differential $d_1:I^*(Fib_{s+1})\to I^*(Fib_s)$ coincides up to sign with the boundary of the normalized cochain complex $N^\bullet(I^*Y^\bullet)$. Hence, the $I^*$-spectral sequence collapses with vanishing $E_2$-term except for $I^*Y=\stabhom{Y}{I}_*$ in degree $0$.

There is also a spectral sequence belonging to the associated tower $\{\widetilde{Y}^s\}$ over $Y$ (see Def. \ref{associated}), that has the same collapsing $E_2$-term. 
By induction on $s$ we can prove with a diagram chase using our knowledge of the $E_2$-term, that the maps $\stabhom{\widetilde{Y}^s}{I}\to\stabhom{\widetilde{Y}^{s+1}}{I}$ are zero for $I\aus\mc{G}{\rm -inj}$. This implies $\colim_s\stabhom{\widetilde{Y}^s}{I}=0$ for all $I\aus\mc{G}{\rm -inj}$. 

To conclude the vanishing of $\colim_s\stabhom{\widetilde{Y}^s}{B}$ for all $B\aus\mc{B}(\mc{G}{\rm -inj})$, we filter the thick subcategory by giving objects in \mc{G}-inj filtration index $0$, and in each succesive step we adjoin retracts and objects built from distinguished triangles. Since colimits behave well we can prove our assertion now by induction on the filtration index.

Property $(ii)$ follows now from the exact sequence
\begin{equation*}
    ...\to\colim_s\stabhom{\widetilde{Y}^s}{B}_{*+1}\to\colim_s\stabhom{\Tot_sY^\bullet}{B}_*\to\stabhom{Y}{B}_*\to\colim_s\stabhom{\widetilde{Y}^s}{B}_*\to ...
\end{equation*}
that is derived from diagram \Ref{überunter}.
\end{beweis}

\section{Nilpotent and injective completion}

We repeat the definition of a $F$-nilpotent resolution and of the $F$-nilpotent completion from \cite{Bou:loc} in the setting of the stable homotopy category of spectra \ho{\mc{S}}. 

\begin{Def}  

Let $F$ be a ring spectrum, or in other words a monoid in the monogenic stable homotopy category \ho{\mc{S}}. We consider the class of spectra of the form $\sm{F}{X}$, where $X$ is some spectrum. We call the objects of the thick subcategory generated by this class {\bf\boldmath $F$-nilpotent spectra\unboldmath}, the subcategory itself we denote by \mc{N}.
\end{Def}

\begin{Def}  

We call a tower $\{W_s\}$ under a spectrum $W$ an {\bf\boldmath $F$-nilpotent resolution of $W$\unboldmath}, if $\{W_s\}$ is \mc{N}-resolution in the sense of definition \Ref{B-Aufl}. 
Bousfield in \cite{Bou:loc} calls 
\begin{align}
      \holim_s W_s =: F^{\wedge}W \notag 
\end{align}
the {\bf $F$-nilpotent completion} of $W$. With our previous conventions it is true that
\begin{center}
      $ F^{\wedge}W = \mc{N}^{\wedge}W $.
\end{center}
We will take over Bousfield's notation.
\end{Def}

Analogously we define an $F$-injective completion. To achieve this we introduce the class of $F$-pseudoinjective objects. Let \mc{D} be a monogenic stable homotopy category.

\begin{Def} \label{Pseudoinjektivit"at}

Let $F$ be a monoid in \mc{D}. The class of {\bf\boldmath $F$-pseudoinjective objects\unboldmath} is the thick subcategory generated by the $F_*$-injective objects, where $F_*=\stabhom{S}{\sm{F}{\frei}}_*$ is the usual homology functor induced by $F$. We abuse language and disregard $F$.
\end{Def}

\begin{Def} 

If we call the class of pseudoinjective objects \mc{P}, then a {\bf\boldmath pseudoinjective resolution\unboldmath} of $W$ in \mc{D} is a \mc{P}-resolution of $W$. We call 
\begin{center}
      $ F^IW := \mc{P}^{\wedge}W $ 
\end{center}
a {\bf\boldmath $F$-injective completion\unboldmath} of $W$.
 
For the {\bf opposite term} (see \ref{dualer Term}) we write:
\begin{center}
    $ W^IF := W^{\wedge}\mc{P} $
\end{center}
\end{Def}

The next lemma states that pseudoinjective resolutions really exist and is proved the same way as \cite[Lemma 5.7]{Bou:loc}.
\begin{lemma} 
Let $Y$ be in \mc{D}, and let $\{Y_s\}$ and $\{Y^s\}$ be $F$-injective resolutions under and over $Y$ \emph{Def. \ref{injektive Auflösung}}. We have: 
\begin{punkt} 
    \item \hspace{2pt} $\{Y_s\}$ is a pseudoinjective resolution.
    \item \hspace{1pt} $\holim\limits_s Y_s \cong F^IY$
    \item $\holim\limits_s Y^s \cong Y^I\!F$
\end{punkt}
\end{lemma} 

\section{Some diagram lemmas}

\begin{Def} \label{linear}
A category \mc{N} is a {\bf simplicial stable model category} if it carries a simplicial model structure, if it is pointed and if for every commutative square in \ho{\mc{N}} the fact, that it is homotopy cartesian, is equivalent to the fact, that it is homotopy cocartesian.
We say, that {\bf\mc{D} comes from a simplicial stable model category} if there is an \mc{N} with $\mc{D}\cong\ho{\mc{N}}$.
Given a stable model category \mc{N} and a finite finite-dimensional category \mc{C} the functor category $\mc{N}^{\mc{C}}$ can be given a stable model structure. 
We can form the family $(\mc{K}_C)$ with 
\begin{center}
      $\mc{K}_C := \ho{\mc{N}^C}$. 
\end{center}
\end{Def}

\begin{Bem}
The families $\mc{K}_C$ considered in \ref{linear} are special cases of systems of triangulated diagram categories developed in \cite{Fra:uni}.
\end{Bem}

\begin{lemma} \label{lange Sub-Sequenz}

There is a long exact sequence
\begin{align}
          & & & & & & \emph{\gradhom{\mc{K}_{\underline{0}}}{q}{A_0}{B_0}} & \nonumber \\ 
    \ldots & \to & \emph{\gradhom{\mc{K}_{\underline{0}}}{q-1}{A_0}{B_1}} & \to & \emph{\gradhom{\mc{K}_{\underline{1}}}{q}{A}{B}} & \to & \oplus \hspace{1.2cm} & \to \emph{\gradhom{\mc{K}_{\underline{0}}}{q}{A_0}{B_1}} \to \ldots \nonumber \\
          & & & & & &  \emph{\gradhom{\mc{K}_{\underline{0}}}{q}{A_1}{B_1}} .& \nonumber  
\end{align}
\end{lemma}

\begin{beweis}
According to \cite{DK:function-complexes}, where we can substitute the simplicial sets used there by our simplicial stable model category \mc{N}, there is a homotopy pull back square
\begin{center}
    $\xymatrix{ \text{map}_{\mc{N}^{\underline{1}}}(A,B) \ar[r]\ar[d] & \text{map}_{\mc{N}}(A_0,B_0) \ar[d] \\
                \text{map}_{\mc{N}}(A_1,B_1) \ar[r] & \text{map}_{\mc{N}}(A_0,B_1)  }$
\end{center}
when $A$ is cofibrant and $B$ is fibrant. This means, that $A_0$ be cofibrant and $A_0\to A_1$ be a cofibration and $B_1$ be fibrant and $B_0\to B_1$ be a fibration, which can be achieved by replacement arguments. The long exact homotopy sequence from the resulting fiber sequence is the desired one.
\end{beweis}

\begin{Bem} \label{Sub-Spec-Seq}
There is a spectral sequence
\begin{equation*} 
    E_2^{p,q}= \lim_{Sub(C)}\!\!\!\!\mbox{}^p\,\, \gradhom{\mc{K}_{\underline{0}}}{q}{A_i}{B_j}\,\, \Longrightarrow\,\, \gradhom{\mc{K}_{C}}{p+q}{A}{B}, 
\end{equation*}
where $i,j\aus C$ and $A_i$ bzw $B_j$ are the values of the functors $A$ and $B$ at $i$ and $j$ and $i\to j$ is an object in $Sub(C)$, where this is the subdivision category defined as follows:
Objects are the morphisms of $C$, morphisms from $f$ to $g$ are factorizations of $g$ through $f$. 
The lemma \ref{lange Sub-Sequenz} can also be proved by noticing that this spectral sequence collapses for $C=\underline{1}$. This proof is necessary, if we want to include the systems of triangulated diagram categories defined in \cite{Fra:uni}. So lemma \ref{lange Sub-Sequenz} and lemma \ref{kurze Sub-Sequenz}, as well as the rest of this article remain valid, if we are given such a system of triangulated diagram categories. With this proof we can also eliminate the requirement of \mc{N} being simplicial.
\end{Bem}

\begin{lemma} \label{kurze Sub-Sequenz}
Let $\mc{D}=\mc{K}_{\underline{0}}$ come from a simplicial stable model category, and $F_*:\mc{D} \to \mc{A}$ be a homological functor. Consider a morphism $A \to B$ between two objects $A=(A_0 \to A_1)$ and $B=(B_0 \to B_1)$ from $\mc{K}_{\underline{1}}$. If $B_1$ is an $F$-injective object in $\mc{D}=\mc{K}_{\underline{0}}$ and the morphism $F_*A_0\to F_*A_1$ is a monomorphism, then there exist short exact sequences of the form:

\hspace{0.5cm}
$ \begin{array}{cccccc}
         & & & \emph{\gradhom{\mc{K}_{\underline{0}}}{q}{A_0}{B_0}} & &  \\ 
         0 & \to \emph{\gradhom{\mc{K}_{\underline{1}}}{q}{A}{B}} & \to & \oplus  & \to \emph{\gradhom{\mc{K}_{\underline{0}}}{q}{A_0}{B_1}} & \to 0 \\
         & & & \emph{\gradhom{\mc{K}_{\underline{0}}}{q}{A_1}{B_1}} & &  
\end{array} $
\end{lemma}

\begin{beweis}
The long exact sequence in lemma \ref{lange Sub-Sequenz} breaks up in short exact sequences, because the horizontal maps
\begin{center}
    $\xymatrix{ \gradhom{\mc{K}_{\underline{0}}}{q}{A_1}{B_1} \ar[r] \ar[d]_{\cong} & \gradhom{\mc{K}_{\underline{0}}}{q}{A_0}{B_1} \ar[d]^{\cong} \\ 
                \Hom{\mc{A}}{F_{*-q}A_1}{F_*B_1} \ar[r] & \Hom{\mc{A}}{F_{*-q}A_0}{F_*B_1} }$
\end{center} 
are surjective by $F$-injectivity of $B_1$.
\end{beweis}

\begin{kor} \label{funktorielle Liftung}
A morphism $B_0 \to B_1$ in $\mc{D}=\mc{K}_{\underline{0}}$ , where $B_1$ is $F$-injective and the induced morphism $F_*B_0 \to F_*B_1$ is a monomorphism, possesses a lifting to the category $\mc{K}_{\underline{1}}$, which is unique up to unique isomorphism.
\end{kor}


\section{Injective completeness of nilpotent objects}

The conditions in the following theorem were explained in \ref{Def. univ. Koeff.}, \ref{Dualität} and \ref{linear}. Originally the theorem was stated in the language of \cite{Fra:uni} using systems of triangulated diagram categories defined there.

\begin{satz} \label{resultat}  
Let \mc{D} be a monogenic stable homotopy category coming from a stable model category.
Let $F$ be a monoid in \mc{D}, where $F_*F$ is flat as $\pi_*F$-right module and $F$ satisfies the universal coefficient property and the duality property.
Let $Y$ be an $F$-nilpotent object. Then $F$ is injectively complete, which means:
\begin{equation*}
      Y\cong F^IY
\end{equation*}
\end{satz}

\begin{beweis} 
By the distinguished triangle 
\begin{equation*}
    Y^IF\to Y\to F^IY\to\Sigma Y^IF
\end{equation*}
from \ref{Dualität der Vervollst.} the isomorphism is equivalent to the statement of the following lemma \ref{Behauptung 2}. 
\end{beweis}

\begin{lemma} \label{Behauptung 2}
If $Y$ is nilpotent then $Y^IF \cong 0$.
\end{lemma}

We have to derive some auxiliary lemmas first.
 
\begin{lemma} \label{Kofaser}
Let $M$ be an $F$-module and let $F_*M\to F_*F\otimes_{\pi_*F}I$ be a monomorphism to an injective $F_*F$-comodule like in cor. \emph{\ref{Tensor-Auflösung}}. Let $E_I$ be the corresponding Eilenberg-MacLane-object, and let 
\begin{equation}
    N \to M \to E_I \to \Sig N  \nonumber
\end{equation}
be a distinguished triangle. Then there exists an an $F$-module structure on $N$, such that $N\to M$ is a morphism of $F$-modules.
\end{lemma}

\begin{Bem}
These distinguished triangles show up in the construction of an injective resolution like in diagram \Ref{Dreieps^s}. We remind the reader that $E_I$, too, is an $F$-module by \ref{spez-inj-Spek-F-Modul}. 

The construction of the $F$-module structure is based on the following observation: 
The fiber of a morphism in the homotopy category is not a functor which causes the usual problems with such statements, but the fiber of an object in the homotopy category of morphisms (of the underlying stable model category) can be given the structure of a functor. 
The proof consists now in showing that in this very special case the morphism $M\to E_I$ admits a lifting to an object in the homotopy category of morphisms that is in $\mc{K}_{\ul{1}}$ unique up to unique isomorphism in that category using lemma \ref{kurze Sub-Sequenz}. 

Here (lemma \ref{kurze Sub-Sequenz}) we need an underlying stable model category, but any model category \mc{C} giving $\ho{\mc{C}}\cong\mc{D}$ will do. It is also worth noting that we do not assume any kind of fancy structure on our monoid $F$ like $A_\8$ or $E_\8$.
\end{Bem}

\begin{proofof}{\ref{Kofaser}}
The first step is to show that the diagram
\diag{ \sm{F}{M} \ar[r]^{\smknapp{1}{f}} \ar[d]_{\mu_M} & \sm{F}{E_I} \ar[d]^{\mu_I} \\
        M \ar[r]_f & E_I }{Ausgang}
commutes in \mc{D}, where $\mu_M$ and $\mu_I$ are the multiplication maps. 
We apply $F_*$ to \Ref{Ausgang}:
\diag{ F_*\sm{F}{M} \ar[r]^{F_*(\smknapp{1}{f})} \ar[d]_{F_*\mu_M} & F_*\sm{F}{E_I} \ar[d]^{F_*\mu_I} \\
       F_*M \ar[r]_{F_*f} & F_*E_I }{F*Ausgang}
If this diagram commutes, so does \Ref{Ausgang} because $E_I$ is a representing object. Proving commutativity of \Ref{F*Ausgang} is straightforward by writing out the definition of the maps used in the diagram.

In the second step we define a multiplication on the fiber. Let \mc{S} be a stable model for \mc{D}.
We consider $\smknapp{1}{f}$ and $f$ as objects in $\mc{D}^{\underline{1}}$, the category of morphisms in \mc{D}, and $(\mu_M,\mu_I)$ as a morphism between these objects. 
The objects $\smknapp{1}{f}$ and $f$ lift to objects in the category $\mc{K}_{\underline{1}}=\ho{\mc{S}^{\ul{1}}}$ (defined in \ref{linear}), and lemma \ref{kurze Sub-Sequenz} applies to these liftings. It follows that there is a unique morphism in $\mc{K}_{\underline{1}}$, which induces the morphism $(\mu_M,\mu_I)$ between the chosen liftings of $f$ and $\sm{1}{f}$, i.e. the commutative square \Ref{Ausgang}.

We pick one lifting of $f$ and call it $\tilde{f}$, we can think of it as represented by a morphism in \mc{S}. In fact it is true, that the lifting $\tilde{f}$ itself is in $\mc{K}_{\ul{1}}$ unique up to unique isomorphism, as we have seen in cor. \ref{funktorielle Liftung}. But this is not needed here.

We call $(\mathfrak{Hofi} \tilde{f})_0$, the vertex $0$ of $\mathfrak{Hofi} \tilde{f}\in \mc{K}_{\underline{2}}$, $N$. As an object in $\mc{K}_{\underline{0}}$ it is a fiber of $f$, i.e. there is a distinguished triangle
\begin{center}
    $ N \to M \toh{f} E_I \to \Sig N $.
\end{center}

But then $\smknapp{1}{\tilde{f}}$ is a lifting of $\smknapp{1}{f}$ with fiber $(\mathfrak{Hofi} \smknapp{1}{\tilde{f}})_0=\sm{F}{N}$. We call the uniquely determined lifting of $(\mu_M,\mu_I)$ $\phi$. Functorially we can form $\mathfrak{Hofi} \phi \in \Hom{\mc{K}_{\underline{2}}}{\mathfrak{Hofi}\hrl\smknapp{1}{\tilde{f}}\hrr}{\mathfrak{Hofi}\hrl\tilde{f}\hrr}$. This is a diagram:
\begin{center}
    $\xymatrix{ \sm{F}{N} \ar[r] \ar[d]_{(\mathfrak{Hofi}\phi)_0} & \sm{F}{M} \ar[r]^{\smknapp{1}{\tilde{f}}} \ar[d]^{(\mathfrak{Hofi}\phi)_1}
      & \sm{F}{E_I} \ar[d]^{(\mathfrak{Hofi}\phi)_2} &  \mathfrak{Hofi}(\smknapp{1}{\tilde{f}}) \ar[d]^{\mathfrak{Hofi} \phi} \\
       N \ar[r] & M \ar[r]_{\tilde{f}} & E_I &  \mathfrak{Hofi}(\tilde{f})  }$
\end{center}
$(\mathfrak{hofi}\phi)_0$ is a morphism in $\mc{K}_{\underline{0}}=\mc{D}$. We define the multiplication $\mu_N$:$\sm{F}{N} \to N$ just as $(\mathfrak{hofi}\phi)_0$. Define the unit $\eta_N:S\to N$ analogously. 

In the last step we have to convince ourselves that $\mu_N$ is associative and possesses a unit. The diagram
\diag{ M \ar[r]^f \ar@{=}[d] & E_I \ar@{=}[d] \\
       M \ar[r]^f  & E_I  }{f-Identität}
possesses by \ref{kurze Sub-Sequenz} a unique lifting. Consider the diagram:
\diag{ M \ar[r] \ar[d]_{\smknapp{\eta}{1}} & E_I \ar[d]^{\smknapp{\eta}{1}} \\
       \sm{F}{M} \ar[r] \ar[d]_{\mu_M} & \sm{F}{E_I} \ar[d]^{\mu_I} \\
       M \ar[r] & E_I  }{philam}
Let $\lam$ be an arbitrary lifting of $(\smknapp{\eta}{1},\smknapp{\eta}{1})$, so that $\phi\lam$ is a lifting of \Ref{philam}. By \ref{kurze Sub-Sequenz} it is unique. 
Then diagram \Ref{f-Identität} states that $id_f$ equals $\phi\lam$ as morphisms in $\mc{K}_{\underline{1}}$, hence by applying the functor $\hofi(\frei )_0$ we arrive at the following equalities of morphisms in $\mc{K}_{\underline{0}}$:
\begin{center}
   $\begin{array}{cll}
       id_N &= (\hofi(id_f))_0 &= (\hofi(\phi\lam))_0  \\
            &= (\hofi(\phi)\, \hofi(\lam))_0 &= (\hofi(\phi))_0\, (\hofi(\lam))_0 \\
            &= \mu_N (\sm{\eta_N}{1}) &
    \end{array}$
\end{center}
$N \stackrel{\smknapp{\eta}{1}}{\longrightarrow} \sm{F}{N}$ is a representative of $(\hofi(\lam))_0$, and it follows:
\begin{equation}
    \mu_N (\smknapp{\eta}{1}) = id_N \nonumber
\end{equation}

To prove associativity we observe that the following diagram commutes in \mc{D}:
\begin{center}
    $\xymatrix{ \sm{F}{\sm{F}{M}} \ar[r] \ar[d]_{\psi_0} & \sm{F}{\sm{F}{E_I}} \ar[d]^{\psi_1} \\
       M \ar[r] & E_I  }$
\end{center}
where
\begin{equation}
      \psi_0:=\mu_M(\smknapp{\mu_F}{1})=\mu_M(\smknapp{1}{\mu_M}) \nonumber
\end{equation}
and
\begin{equation}
      \psi_1:=\mu_I(\smknapp{\mu_F}{1})=\mu_I(\smknapp{1}{\mu_I}) . \nonumber
\end{equation}
Of course we have used the associativity of $\mu_M$ and $\mu_I$. Again the lifting is unique because of lemma \ref{kurze Sub-Sequenz}, and we call it $\psi$. Let $\kap^1$ be a lifting of $(\smknapp{\mu_F}{1},\smknapp{\mu_F}{1})$ and let $\kap^2$ be a lifting $(\smknapp{1}{\mu_M},\smknapp{1}{\mu_I})$. $\psi$ factors in the following way:
\begin{equation}
    \phi \kap^1 = \psi = \phi \kap^2 \nonumber
\end{equation} 
We form $(\hofi(\frei))_0$ and obtain:
\begin{equation}
    \mu_N (\kap^1)_0 = \mu_N (\kap^2)_0 \nonumber
\end{equation}
We can choose $\smknapp{\mu_F}{1}:\sm{F}{\sm{F}{N}}\to \sm{F}{N}$ and $\smknapp{1}{\mu_N}:\sm{F}{\sm{F}{N}} \to \sm{F}{N}$ as representatives of $(\kap^1)_0$ and $(\kap^2)_0$, whence we have proved associativity.

That the morphisms are morphisms of $F$-modules is clear from the construction.\end{proofof}

\begin{kor} \label{Modulauflösung}
Let $Y = \sm{F}{X}$ for some $X$ in \mc{D}. Then there exists an injective resolution over $Y$ by $F$-modules. 
\end{kor}

\begin{beweis}
$Y^0\!=\!Y\!=\!\sm{F}{X}$ is an $F$-module. Corollary \ref{Tensor-Auflösung} and lemma \ref{Kofaser} provide the inductive step.
\end{beweis}

\begin{proofof}{\ref{Behauptung 2}}
The thick subcategory \mc{N} of $F$-nilpotent objects has a filtration. We prove the claim by induction over the filtration index.

Assume $Y$ to be an object of $\mc{N}_0$, so there is $X$ in \mc{D} with $Y \cong \sm{F}{X}$. Then we know by corollary \ref{Modulauflösung}, that there is an injective resolution $\kl Y^s\kr$ over $Y$ with all $Y^s$ having an $F$-module structure. By corollary \ref{univ. Koeff.} and definition \ref{m} we arrive at a commutative diagram of the form
\begin{center}
    $\xymatrix{ F_*Y^{s+1} \ar[r]^0 \ar[d]^{\cong} & F_*Y^{s} \ar[d]^{\cong} \\
       F_*F \otimes \pi_*Y^{s+1} \ar[r] & F_*F \otimes \pi_*Y^s  }$
\end{center} 
where the horizontal maps vanish by construction of an injective resolution. Because $F_*F$ is faithfully flat, it follows, that actually the maps
\begin{equation}
    \pi_*Y^{s+1} \toh{0} \pi_*Y^s  \label{0}
\end{equation}
are zero. $\pi_*$ has a Milnor sequence
\begin{equation}
    0 \to \Rlim_{s} \pi_{*+1}Y^s \to \pi_*\holim_{s} Y^s \to \lim_{s} \pi_*Y^s \to 0\, . \nonumber
\end{equation}
The outer terms vanish as we can see in \Ref{0}, therefore $\pi_* \hrl\holim Y^s\hrr = \pi_* \hrl Y^IF\hrr$ vanishes. We get
\begin{center}
    $Y^IF \cong \holim\limits_s Y^s \cong 0$
\end{center}
The inductive step is proved by using Theorem \ref{Eigenschaften} and contains no further difficulties.
\end{proofof}

\section{Comparison of the two completions}

In this last section let \mc{G} be {\bf the class of $F$-injective objects}, and let $\mc{G}'$ be {\bf the class of all $F$-modules}.

\begin{lemma} \label{inj folgt nil}
Every $F$-injective object $I$ is a retract of $\sm{F}{I}$, so every $F$-injective object is $F$-nilpotent.
\end{lemma}

\begin{kor} \label{missinglink}
Every \mc{G}-injective object is $\mc{G}'$-injective, and every $\mc{G}'$-injective object is \mc{G}-complete.
\end{kor}

\begin{beweis}
$(i)$ follows from \ref{inj folgt nil}, and $(ii)$ follows from \ref{resultat}.
\end{beweis}

\begin{lemma} \label{G}
For every $X$ in \mc{D} there are natural equivalences
\begin{eqnarray*}
      \hat{L}_{\mc{G}}X\cong F^IX & \text{and} & \hat{L}_{\mc{G'}}X\cong F^{\wedge}X.
\end{eqnarray*}
\end{lemma}

\begin{beweis}
This is just a specialization of lemma \ref{Tot}.
\end{beweis}

\begin{satz} \label{inj=nil}
Let $F$ be a monoid as in \emph{\ref{resultat}}. We assume, that the underlying model category is simplicial stable and left proper. Then for every object $X$ there is a natural equivalence
\begin{equation*}
    F^IX\cong F^{\wedge}X
\end{equation*} 
\end{satz}

\begin{beweis}
Take \ref{G} and \ref{missinglink} together with theorem \ref{bousfield}.
\end{beweis}

\section{Application}

The only rather immediate applications I can prove are the following two lemmas. The conditions in the following lemma were explained in \ref{Def. univ. Koeff.}, \ref{Dualität} and \ref{linear}.
\begin{lemma} \label{Konvergenz}
Let $F$ be a monoid with $F_*F$ flat over $\pi_*F$ satisfying the universal coefficient and the duality property. Let the underlying model category be simplicial, stable and left proper.
Then the modified Adams spectral sequence $E_r^{*,*}(X,Y)$ converges strongly to $\emph{\stabhom{X}{F^{\wedge}Y}}$ if and only if $\Rlim\limits_r E_r^{*,*}(X,Y)=0$.
\end{lemma}

\begin{beweis}
By the arguments in \cite{Boa:ccss} or in \cite[p. 263]{BK:lim} it follows that the result is true for $\stabhom{X}{F^IY}$, and $F^IY\cong F^{\wedge}Y$ by \ref{inj=nil}.
\end{beweis} 

\begin{Bem}
It follows in particular, if the injective dimension of the target category is finite, then the modified Adams spectral sequence converges strongly. We also have to following result about the classical version of the Adams spectral sequence for a monoid $F$. 
\end{Bem}

\begin{kor} \label{klassische Konvergenz}
Let $F$ be as above and suppose, that $F_*X$ is projective as $\pi_*F$-module. If the injective dimension of the category of $F_*F$-comodules is finite, then the Adams spectral sequence $E_r^{*,*}(X,Y)$ converges strongly to $\emph{\stabhom{X}{F^{\wedge}Y}}$.
\end{kor}

\begin{beweis}
There is a \mc{G}-equivalence $\sm{E}{\sm{\ol{E}^\bullet}{Y}}\to E\hrl I^\bullet\hrr$ for \mc{G} the class of $F$-injectives, where the first cosimplicial object is used to construct the Adams spectral sequence, while the second one is used to construct the modified version. For projective $F_*X$ the functor $\stabhom{X}{\frei}$ maps this equivalence to an equivalence of cosimplicial objects in the abelian category of $F_*F$-comodules. Thus the resulting spectral sequences are isomorphic from the $E_2$-term on. This is also proved in \cite[prop. 1.9]{Dev:brown-comenetz}.
The result now follows from \ref{Konvergenz}.
\end{beweis}

As an example we can generalize \cite[5.3.]{Hov-Sad:invertible} where rather difficult results were used to show $E(n)$-prenilpotency. Here $E(n)$ denotes the $n$-th Johnson-Wilson spectrum.
\begin{kor} 
Let $p$ be a prime with $p\uber n+1$. Let $X$ and $Y$ be $E(n)$-local spectra, where $E(n)_*X$ is projective as $E(n)_*$-module.
Then $E(n)^{\wedge}Y\cong Y$ and the $E(n)$-based Adams spectral sequence $E_r^{*,*}(X,Y)$ converges strongly to $\emph{\stabhom{X}{Y}}$.
\end{kor}

\begin{beweis}
For $p\uber n+1$ it is shown in \cite[5.1.]{Hov-Sad:invertible} that the injective dimension of $E(n)_*E(n)$-comodules is $n^2+n$, in particular finite. Then the first claim follows from \ref{inj=nil} and the second from \ref{klassische Konvergenz}.
\end{beweis}

\begin{Def}
We will denote the {\bf Bousfield localization} functor by $\mc{L}_F(\frei)$.
An object in $\mc{D}$ is called {\bf preinjective}, if its Bousfield-localization is pseudoinjective (Def. \ref{Pseudoinjektivit"at}) respectively.
\end{Def}

\begin{Bem} 
The relation of "preinjective" to "pseudoinjective" is as "prenilpotent" to "nilpotent".
\end{Bem}

\begin{lemma} 
Assume additionally, that $Y$ is preinjective. Then the modified Adams spectral sequence converges strongly to $\emph{\stabhom{X}{\mc{L}_FY}}$. There is $s_0 \aus \mathbb{N}$ and $\varphi: \mathbb{N} \to \mathbb{N}$, such that for every $X$ $E_{\8}^s(X,Y) = 0$, if $s\uber s_0$, and $E_r^s(X,Y) = E_{\8}^s(X,Y)$, if $r\uber\varphi(s)$.
\end{lemma}

\begin{beweis}
The constant tower $\{\mc{L}_FY\}$ is a pseudoinjective resolution. This implies, that $\kl\stabhom{X}{Y_s}\kr$ is pro-isomorphic to a constant tower. 
\end{beweis}

\begin{Bem} 
The assertion of the last lemma is also true, if we assume that $F_*X$ is projective over $\pi_*F$ and $Y$ is prenilpotent. This yields \cite[thm. 6.10]{Bou:loc}.
\end{Bem}

\end{document}